\documentclass[11pt]{article}
\usepackage{amsthm,amstext}
\usepackage{amsmath}
\usepackage{graphicx}
\usepackage{amsfonts}
\usepackage{amssymb}
\theoremstyle{plain}
\newtheorem{theorem}{Theorem}[section]
\newtheorem{lemma}[theorem]{Lemma}
\newtheorem{proposition}[theorem]{Proposition}

\theoremstyle{definition}
\newtheorem{definition}[theorem]{Definition}
\newtheorem{corollary}[theorem]{Corollary}
\newtheorem{example}{\sc Example}
\theoremstyle{remark}
\newtheorem{remark}{\sc Remark}

\date{}
\title{\bf On Fuzzy Bi-ideals and Fuzzy Quasi Ideals \\in $\Gamma $-Semigroups}\vspace{.25 in}
\author{ \textbf{\footnote{The Research is partially supported by DST(Department of Science and Technology, Govt. of India)-PURSE,} $^{1}$Sujit Kumar Sardar,} \textbf{$^{2}$Samit Kumar Majumder}\\\\
and\\\\ \textbf{$^{3}$Soumitra Kayal} \\\\
$^{1,3}$Department of Mathematics, Jadavpur\\
University, Kolkata-700032, INDIA\\
$^2$Tarangapur N.K High School, Tarangapur,\\
Uttar Dinajpur, West Bengal-733129, INDIA\\
{\tt  $^1$sksardarjumath@gmail.com}\\
{\tt $^2$samitfuzzy@gmail.com}\\
{\tt $^3$soumitrakayal.ju@gmail.com}
 }
\begin{document}
\maketitle

\begin{abstract}
The purpose of this paper is to investigate some properties of fuzzy ideals and fuzzy bi-ideals in $\Gamma $-semigroups and to introduce the notion of fuzzy quasi ideals in $\Gamma $-semigroups. Here we also characterize a regular $\Gamma$-semigroup in terms of fuzzy quasi ideals.


\end{abstract}

\section{Introduction}
A semigroup is an algebraic structure consisting of a non-empty set $S$ together with an associative binary operation\cite{H}. The formal study of semigroups began in the early 20th century. Semigroups are important in many areas of mathematics, for example, coding and language theory, automata theory, combinatorics and mathematical analysis. Sen and Saha in \cite{SS}, defined $\Gamma$-semigroups as a generalization of semigroups. $\Gamma$-semigroups have been analyzed by a lot of mathematicians, for instance by Chattopadhyay\cite{Ch1, Ch2}, Dutta and Adhikari\cite{D1,D2}, Hila\cite{H1, H2}, Chinram\cite{Chin,Chin1}, Saha\cite{S}, Sen et al.\cite{SC, SSa,SSe}, Seth\cite{Se}.\\
\indent After the introduction of fuzzy sets by  Zadeh\cite{Z}, reconsideration of the concept of classical mathematics began. On the other hand, because of the importance of group theory in mathematics, as well as its many areas of application, the notion of fuzzy subgroups was defined by Rosenfeld\cite{R} and its structure was investigated. Das characterized fuzzy subgroups by their level subgroups in \cite{Das}. \\
\indent Hong, Jun and Meng\cite{Jun} considered the fuzzification of interior ideals in semigroups. Sardar and Majumder\cite{S2,S4} characterized interior ideals(along with B. Davvaz\cite{S5}), ideals, prime(along with D. Mandal\cite{S3}) and semiprime ideals, ideal extensions(along with T.K. Dutta\cite{D3,D4}) of a $\Gamma$-semigroup in terms of fuzzy subsets. This paper is a sequel to our study \cite{D3,D4,S2,S3,S4,S5} of fuzzification of $\Gamma $-semigroups. In this paper we have investigated some properties of fuzzy ideals, fuzzy bi-ideals and fuzzy $(1,2)$-ideals and characterized a $\Gamma $-semigroup which is left$($right$)$ simple, left$($ right$)$ duo, left$($right$)$ regular, intra-regular, regular in terms of fuzzy ideals and fuzzy bi-ideals.\\
\indent The concept of quasi-ideals in semigroups was introduced
in $1956$ by O. Steinfeld\cite{St}. The class of quasi-ideals in semigroups is a generalization of one sided ideals in semigroups. It is well known that the intersection of a left ideal and a right ideal of a semigroup $S$ is a quasi ideal of $S$ and every quasi ideal of $S$ can be obtained in this way. Nobuaki Kuroki\cite{K} characterized completely regular semigroup and semigroup in terms of fuzzy semiprime quasi ideals. He also studied the properties of fuzzy bi-ideals in semigroups. We have introduced the concept of fuzzy quasi ideals in $\Gamma $-semigroups and studied some
of its important properties in this paper. Lastly we have obtained a characterization theorem for regular $\Gamma$ semigroups in terms of fuzzy quasi ideals.
\section{Preliminaries}

In this section we discuss some elementary definitions that we use in the sequel.

Let $S=\{x,y,z,\ldots\}$ and $\Gamma=\{\alpha,\beta,\gamma,\ldots\}$ be two non-empty sets. Then $S$ is called a {\it $\Gamma$-semigroup}\cite{SS} if there exists a mapping $S \times\Gamma\times S $  $\rightarrow S$ (images to be denoted by $ a\alpha b$) satisfying
\begin{itemize}
\item[(1)] $x\gamma y \in S$,
\item[(2)] $(x\beta y)\gamma z=x\beta (y\gamma z),\forall x,y,z\in S,\forall\gamma\in\Gamma.$
\end{itemize}

\begin{remark}
The $\Gamma$-semigroup introduced by Sen and Saha\cite{SS} may be called {\it one sided $\Gamma$-semigroup}. Later Dutta and Adhikari\cite{D1} introduced both sided $\Gamma$-semigroup where the operation $\Gamma \times S \times \Gamma$ to $\Gamma$ was also taken into consideration. They defined operator semigroups for such $\Gamma$-semigroups.
\end{remark}

A non-empty subset $I$ of a $\Gamma$-semigroup $S$ is said to be a {\it subsemigroup} of $S$ if $I \Gamma I \subseteq I$. A subsemigroup $I$ of a $\Gamma $-semigroup $S$ is called a {\it bi-ideal} of $S$ if $I\Gamma S\Gamma I\subseteq I.$ A subsemigroup $I$ of a $\Gamma $-semigroup $S$ is called a $(1,2)$-{\it ideal} of $S$ if $I\Gamma S\Gamma I\Gamma I\subseteq I.$

Let $S$ be a $\Gamma $-semigroup. A non-empty subset $Q$ of $S$ is called {\it quasi ideal}\cite{Chin} of $S$ if $S\Gamma Q\cap Q\Gamma S\subseteq Q.$ Let $Q$ be a quasi ideal of a $\Gamma$-semigroup $S.$ Then $Q\Gamma Q\subseteq S\Gamma Q\cap Q\Gamma S\subseteq Q.$ Hence $Q$ is a subsemigroup of $S.$

\begin{example}
Let $\Gamma =\{5,7\}$. For any $x,y\in \mathbb {N}$ and $\gamma \in \Gamma$, define $x\gamma y=x\cdot \gamma \cdot y$ where $\cdot$ is the usual multiplication on $\mathbb{N}$. Then $\mathbb{N}$ is a $\Gamma$-semigroup.
\end{example}

\begin{example}
Let $S=[0,1]$ and $\Gamma=\{\frac{1}{n} \ | \ n \ {\rm is \ a \ positive \ integer} \}$. Then $S$ is a $\Gamma$-semigroup under the usual multiplication. Next, let $K=[0,1]$. We have that $K$ is a non-empty subset of $S$ and $a \gamma b \in K$ for all $a,b\in K$ and $\gamma \in \Gamma$. Then $K$ is a sub $\Gamma$-semigroup of $S$.
\end{example}

\begin{example}
\cite{Chin1} Let $S$ be a semigroup and $\Gamma =\{1\}.$ Let us define a mapping $S\times\Gamma\times S\rightarrow\Gamma$ by $a1b= ab\forall a,b\in S.$ Then $S$ is a $\Gamma$-semigroup. Let $B$ be a bi-ideal of the semigroup $S.$ Thus $BSB\subseteq B.$ Since $\Gamma= \{1\},B\Gamma S\Gamma B = BSB\subseteq B.$ Hence $B$ is a
bi-ideal of the $\Gamma$-semigroup $S.$
\end{example}

\begin{example}
\cite{Chin} Let $S$ be a semigroup and $\Gamma$ be any non-empty set. Let us define a mapping $S\times\Gamma\times
S\rightarrow\Gamma$ by $a\gamma b=ab\forall a,b\in S,\forall\gamma\in\Gamma.$ Then $S$ is a $\Gamma$-semigroup. Let
$Q$ be a quasi ideal of the semigroup $S.$ Then $SQ\cap QS\subseteq Q.$ Again we see that $S\Gamma Q\cap Q\Gamma S=SQ\cap QS\subseteq Q.$ Hence $Q$ is a quasi ideal of the $\Gamma$-semigroup $S.$
\end{example}

Throughout the paper unless otherwise stated $S$ will denote a one-sided $\Gamma$-semigroup.\\

A {\it left} ({\it right}) {\it ideal}\cite{SS} of a $\Gamma$-semigroup $S$ is a non-empty subset $I$ of $S$ such that
$S\Gamma I \subseteq I$ ($I \Gamma S \subseteq I$). If $I$ is both a left and a right ideal of a $\Gamma$-semigroup $S$, then we say that $I$ is an {\it ideal} of $S$.

Let $S$ and $T$ be two $\Gamma$-semigroups. A mapping $f:S\rightarrow T$ is called a {\it homomorphism} if
$$f(x\gamma y)=f(x)\gamma f(y)\ {\rm for \ all} \ x,y\in S \ {\rm and} \ \gamma\in\Gamma.$$

Though this is defined for both sided $\Gamma$-semigroup in \cite{D1}, the same definition can be adopted for one-sided $\Gamma$-semigroup of Sen and Saha\cite{SS}.\\

A {\it fuzzy subset}\cite{Z} in $S$ is a function $\mu : S\longrightarrow [0,1]$.

Let $\mu $ and $\sigma $ be two fuzzy subsets of a $\Gamma$-semigroup $S.$ Then we define\cite{S2} the following:

\begin{itemize}
\item[(1)] $(\mu \cap \sigma )(x)=\min \{\mu (x),\sigma (x)\}\forall x\in S,$
\item[(2)] $(\mu \cup \sigma )(x)=\max \{\mu (x),\sigma (x)\}\forall x\in S,$
\item[(3)] $(\mu \circ \sigma )(x)=\left\{
\begin{array}{c}
\underset{x=y\gamma z}{\sup }[\min \{\mu (y),\sigma (z)\}] \\
0,\text{ \ \ \ \ Otherwise}%
\end{array}%
\right.. $
\end{itemize}

A non-empty fuzzy subset $\mu$ of a $\Gamma$-semigroup $S$ is called a {\it fuzzy left ideal}\cite{S2} of $S$ if $\mu(x\gamma y)\geq\mu(y)$ $\forall x,y\in S$, $\forall\gamma\in\Gamma.$

A non-empty fuzzy subset $\mu$ of a $\Gamma$-semigroup $S$ is called a {\it fuzzy right ideal}\cite{S2} of $S$ if $\mu(x\gamma y)\geq\mu(x)\forall x,y\in S, \forall\gamma\in\Gamma.$

A non-empty fuzzy subset $\mu$ of a $\Gamma$-semigroup $S$ is called a {\it fuzzy ideal}\cite{S2} of $S$ if it is both a fuzzy left ideal and a fuzzy right ideal of $S.$

In a $\Gamma $-semigroup $S$ following are equivalent: $(1)$ $\mu$ is a fuzzy left$($right$)$ ideal of $S,$ $(2)$ $\chi\circ\mu\subseteq\mu(\mu\circ\chi\subseteq\mu)$(\cite{S2}).

\section{Fuzzy Subsemigroup and Fuzzy Bi-ideal}

\begin{definition}
A non-empty fuzzy subset $\mu $ of a $\Gamma$-semigroup $S$ is called a fuzzy subsemigroup of $S$ if $\mu(x\gamma y)\geq\min\{\mu(x),\mu(y)\}$ $\forall x,y\in S$, $\forall\gamma\in\Gamma.$
\end{definition}

\begin{definition}
A fuzzy subsemigroup $\mu $ of a $\Gamma$-semigroup $S$ is called a fuzzy bi-ideal of $S$ if $\mu(x\beta s\gamma y)\geq\min\{\mu(x),\mu(y)\}\forall x,s,y\in S,\forall\beta,\gamma\in\Gamma.$
\end{definition}

\begin{definition}
A fuzzy subsemigroup $\mu $ of a $\Gamma$-semigroup $S$ is called a fuzzy $(1,2)$-ideal of $S$ if $\mu(x\alpha\omega\beta (y\gamma z))\geq\min\{\mu(x),\mu(y),\mu(z)\}$ $\forall x,\omega ,y,z\in S,$ $\forall\alpha,\beta,\gamma\in\Gamma.$
\end{definition}

\begin{theorem}
Let $I$ be a non-empty subset of a $\Gamma $-semigroup $S$ and $\lambda _{I}$ be the characteristic function of $I.$ Then $I$ is a $\Gamma $-subsemigroup of $S$ if and only if $\lambda _{I}$ is a fuzzy subsemigroup of $S$.
\end{theorem}

\begin{proof}
Let $I$ be a subsemigroup of $S$. Let $x,y\in S$ and $\gamma\in\Gamma .$ Then $x\gamma y\in I$ if $x,y\in I.$ It follows that $\lambda_{I}(x\gamma y)=1=\lambda_{I}(x)=\lambda_{I}(y)=\min\{\lambda_{I}(x),\lambda_{I}(y)\}.$ Let either $x\notin I$ or $y\notin I$. Then Case-$(1):$ If $x\gamma y\notin I,$ then $\lambda_{I}(x\gamma y)= 0=\min\{\lambda_{I}(x),\lambda_{I}(y)\}.$ Case-$(2):$ If $x\gamma y\in I$, then $\lambda_{I}(x\gamma y)=1\geq 0=\min \{\lambda _{I}(x),\lambda_{I}(y)\}.$ Thus $\lambda_{I}$ is a fuzzy subsemigroup of $S$.

Conversely, let $\lambda_{I}$ be a fuzzy subsemigroup of $S.$ Let $x,y\in I,$ then $\lambda_{I}(x)=\lambda_{I}(y)=1.$ Thus $\lambda_{I}(x\gamma y)\geq\min\{\lambda_{I}(x),\lambda_{I}(y)\}=1$ $\forall\gamma\in\Gamma.$ Hence $x\gamma y\in I$ $\forall\gamma\in\Gamma.$ So $I$ is a subsemigroup of $S$.
\end{proof}

\begin{theorem}
Let $I$ be a non-empty subset of a $\Gamma$-semigroup $S$ and $\lambda_{I}$ be the characteristic function of $I.$ Then $I$ is a bi-ideal of $S$ if and only if $\lambda_{I}$ is a fuzzy bi-ideal of $S$.
\end{theorem}

\begin{proof}
Let $I$ be a bi-ideal of a $\Gamma$-semigroup $S$ and $\lambda_{I}$ be the characteristic function of $I.$ Then $I$ is a subsemigroup of $S.$ Hence by Theorem $3.4,$ $\lambda_{I}$ is a fuzzy subsemigroup of $S.$ Let $x,y,z\in S$ and
$\beta,\gamma\in\Gamma.$ Then $x\beta y\gamma z\in I$ if $x,z\in I.$ It follows that $\lambda_{I}(x\beta y\gamma z)=1=\min\{\lambda_{I}(x),\lambda_{I}(z)\}.$ Let either $x\notin I$ or $z\notin I.$ Then Case-$(1):$ If
$x\beta y\gamma z\notin I,$ then $\lambda_{I}(x\beta y\gamma z)\geq 0=\min\{\lambda_{I}(x),\lambda_{I}(z)\}.$ Case-$(2):$ If $x\beta y\gamma z\in I,$ then $\lambda_{I}(x\beta y\gamma z)=1\geq 0=\min\{\lambda_{I}(x),\lambda_{I}(z)\}.$ Consequently, $\lambda_{I}$ is a fuzzy bi-ideal of $S.$

Conversely, $\lambda_{I}$ is a fuzzy bi-ideal of $S.$ So $\lambda_{I}$ is a fuzzy subsemigroup of $S.$ Hence by Theorem $3.4,$ $I$ is a subsemigroup of $S.$ Let $x,z\in I,$ then $\lambda_{I}(x)=\lambda_{I}(z)=1.$ Thus
$\lambda_{I}(x\beta y\gamma z)\geq\min\{\lambda_{I}(x),\lambda_{I}(z)\}=1$ $\forall y\in S,\forall\beta,\gamma\in\Gamma .$ Hence $x\beta y\gamma z\in I$ $\forall y\in S,\forall\beta,\gamma\in\Gamma.$ So $I$ is a bi-ideal of $S.$
\end{proof}

\begin{theorem}
Let $S$ be a $\Gamma$-semigroup and $\mu$ be a non-empty fuzzy subset of $S.$ Then $\mu$ is a fuzzy subsemigroup of $S$ if and only if $\mu_{t}$ is a subsemigroup of $S$ for all $t\in Im(\mu)$ where $\mu_{t}=\{x\in S:\mu(x)\geq t\}.$
\end{theorem}

\begin{proof}
Let $\mu$ be a fuzzy subsemigroup of $S.$ Let $t\in Im(\mu),$ then there exist some $z\in S$ such that $\mu(z)=t$ and so $z\in\mu_{t}.$ Thus $\mu_{t}\neq\phi .$ Let $x,y\in\mu_{t}.$ Then $\mu(x)\geq t$ and $\mu(y)\geq t$ whence $\min\{\mu(x),\mu(y)\}\geq t.$ Now for $\gamma\in\Gamma,$ $\mu (x\gamma y)\geq\min\{\mu(x),\mu(y)\}.$ Hence $x\gamma y\in\mu _{t}\forall\gamma\in\Gamma,i.e.,\mu_{t}\Gamma\mu_{t}\subseteq \mu_{t}.$ Consequently, $\mu_{t}$ is a
subsemigroup of $S.$

Conversely, suppose $\mu_{t}$'s are subsemigroups of $S,$ for all $t\in Im(\mu).$ Let $x,y\in S$ and $\gamma\in\Gamma.$ Let $\mu(x)=t_{1}$ and $\mu(y)=t_{2}.$ Without any loss of generality suppose $t_{1}\geq t_{2}$. Then $x,y\in\mu_{t_{2}}.$ Then by hypothesis $x\gamma y\in \mu_{t_{2}}.$ Hence $\mu(x\gamma y)\geq t_{2}=\min\{\mu(x),\mu(y)\}.$ Consequently, $\mu$ is a fuzzy subsemigroup of $S.$ Hence the proof.
\end{proof}

\begin{theorem}
Let $S$ be a $\Gamma $-semigroup and $\mu$ be a non-empty fuzzy subset of $S.$ Then $\mu$ is a fuzzy bi-ideal of $S$ if and only if $\mu_{t}$ is a bi-ideal of $S$ for all $t\in Im(\mu)$ where $\mu_{t}=\{x\in S:\mu(x)\geq t\}.$
\end{theorem}

\begin{proof}
Let $\mu$ be a fuzzy bi-ideal of $S.$ Then $\mu$ is a fuzzy subsemigroup of $S.$ Hence by Theorem $3.6,$ $\mu_{t}$ is a subsemigroup of $S,$ $\forall t\in Im(\mu).$ Let $t\in Im(\mu).$ Then there exists some $\alpha\in S$ such that $\mu(\alpha )=t$ and so $\alpha\in\mu_{t}.$ Thus $\mu_{t}\neq\phi .$ Let $x,z\in\mu_{t}.$ Then $\mu(x)\geq t$ and $\mu(z)\geq t$ whence $\min\{\mu(x),\mu(z)\}\geq t.$ Now for $y\in S,\alpha,\beta\in\Gamma,\mu(x\alpha y\beta z)\geq \min\{\mu(x),\mu(z)\}.$ Hence $\mu(x\alpha y\beta z)\geq t$ whence $x\alpha y\beta z\in \mu _{t}.$ Hence $\mu_{t}\Gamma S\Gamma\mu_{t}\subseteq\mu _{t}.$ Consequently, $\mu_{t}$ is a bi-ideal of $S.$

Conversely, suppose $\mu_{t}$'s are bi-ideals of $S,$ for all $t\in Im(\mu).$ Then $\mu_{t}$'s are subsemigroups of $S.$ Hence by Theorem $3.6,$ $\mu$ is a fuzzy subsemigroup of $S.$ Let $x,y,z\in S$ and $\beta,\gamma\in\Gamma.$ Let $\mu(x)=t_{1},\mu(z)=t_{2}.$ Without any loss of generality suppose $t_{1}\geq t_{2}$. Then $x,z\in\mu_{t_{2}}.$ Hence by hypothesis $x\beta y\gamma z\in\mu_{t_{2}}$ whence $\mu(x\beta y\gamma z)\geq t_{2}=\min\{\mu(x),\mu(z)\}.$ Hence $\mu$ is a fuzzy bi-ideal of $S.$ Hence the proof.
\end{proof}

\begin{proposition}
Let $\mu$ and $\lambda$ be two fuzzy subsemigroups$($fuzzy bi-ideals$)$ of a $\Gamma$-semigroup $S.$ Then $\mu\cap\lambda$ is also a fuzzy subsemigroup$($resp. fuzzy bi-ideal$)$ of $S$ provided $\mu\cap\lambda$ is non-empty.
\end{proposition}

\begin{proof}
Let $\mu$ and $\lambda$ be two fuzzy subsemigroups of a $\Gamma$-semigroup $S.$ Let $x,y\in S$ and $\gamma\in\Gamma.$ Then
\begin{align*}
(\mu\cap\lambda)(x\gamma y)&=\min\{\mu(x\gamma y),\lambda(x\gamma y)\}\\
&\geq\min[\min\{\mu(x),\mu(y)\},\min\{\lambda(x),\lambda(y)\}]\\
&=\min[\min\{\mu(x),\mu(y),\lambda(x),\lambda(y)\}]\\
&=\min[\min\{\mu(x),\lambda(x)\},\min\{\mu(y),\lambda(y)\}]\\
&=\min\{(\mu\cap\lambda)(x),(\mu\cap\lambda)(y)\}.
\end{align*}
Hence $\mu\cap\lambda$ is a fuzzy subsemigroup of $S.$ Similarly we can prove the result for fuzzy bi-ideal.
\end{proof}

\begin{proposition}
Let $f:R\rightarrow S$ be a homomorphism of $\Gamma$-semigroups.

\begin{itemize}
 \item[$(1)$] If $\lambda$ is a fuzzy subsemigroup$($fuzzy bi-ideal$)$ of $S$, then $f^{-1}(\lambda)$ is also a fuzzy subsemigroup$($resp. fuzzy bi-ideal$)$ of $R($where $f^{-1}(\lambda)(r\gamma s):=\lambda(f(r\gamma s))$ for all $r,s\in R$ and $\gamma\in\Gamma),$ provided $f^{-1}(\lambda)$ is non-empty.
 \item[$(2)$] If $f$ is a surjective homomorphism and $\mu$ is a fuzzy subsemigroup$($fuzzy bi-ideal$)$ of $R$ then, $f(\mu)$ is also a fuzzy subsemigroup$($resp. fuzzy bi-ideal$)$ of $S($where $(f(\mu))(x^{^{\prime}}\gamma y^{^{\prime}}):=\underset{f(z)=x^{^{\prime}}\gamma y^{^{\prime}}}{\sup}\mu(z)$ for all $x^{^{\prime}},y^{^{\prime}}\in S$ and $\gamma\in\Gamma).$
\end{itemize}
\end{proposition}

\begin{proof}
$(1)$ Let $\lambda$ be a fuzzy subsemigroup of $S.$ Let $r,s\in R$ and $\gamma\in\Gamma.$ Then

\begin{align*}
(f^{-1}(\lambda))(r\gamma s)&=\lambda(f(r\gamma s))\\
&=\lambda(f(r)\gamma f(s))\\
&\geq\min\{\lambda(f(r)),\lambda(f(s))\}\\
&=\min\{(f^{-1}(\lambda))(r),(f^{-1}(\lambda))(s)\}.
\end{align*}

Hence $(f^{-1}(\lambda))$ is a fuzzy subsemigroup of $R.$\\

$(2)$ Let $\mu$ be a fuzzy subsemigroup of $R.$ Since $(f(\mu))(x^{^{\prime}})=\underset{f(x)=x^{^{\prime}}}{\sup}\mu(x)$ for $x^{^{\prime}}\in S,$ so $f(\mu)$ is non-empty. Let $x^{^{\prime}},y^{^{\prime}}\in S$ and $\gamma\in\Gamma.$ Then

\begin{align*}
(f(\mu))(x^{^{\prime}}\gamma y^{^{\prime}})&=\underset{f(z)=x^{^{\prime}}\gamma y^{^{\prime}}}{\sup}\mu(z)\\
&\geq\underset{%
\begin{array}{c}
f(x)=x^{^{\prime}} \\
f(y)=y^{^{\prime}}%
\end{array}%
}{\sup}\mu(x\gamma y)\\
&\geq\underset{%
\begin{array}{c}
f(x)=x^{^{\prime}} \\
f(y)=y^{^{\prime}}%
\end{array}%
}{\sup}\min\{\mu(x),\mu(y)\}\\
&=\min\{\underset{f(x)=x^{^{\prime}}}{\sup}\mu(x),\underset{f(y)=y^{^{\prime}}}{\sup}\mu(y)\}\\
&=\min\{(f(\mu))(x^{^{\prime}}),(f(\mu))(y^{^{\prime}})\}.
\end{align*}

Consequently, $f(\mu)$ is a fuzzy subsemigroup of $S.$ Similarly we can prove the result for fuzzy bi-ideal.
\end{proof}

\begin{proposition}
Suppose $\theta$ is an endomorphism and $\mu$ is a fuzzy subsemigroup$($fuzzy bi-ideal$)$ of a $\Gamma$-semigroup $S.$ Then $\mu\lbrack\theta ]$ is also a fuzzy subsemigroup$($resp. fuzzy bi-ideal$)$ of $S,$ where $\mu[\theta](x):=\mu(\theta(x))\forall x\in S.$
\end{proposition}

\begin{proof}
Let $x,y\in S$ and $\gamma\in\Gamma .$ Then
\begin{align*}
\mu\lbrack\theta](x\gamma y)&=\mu(\theta(x\gamma y))\\
&=\mu(\theta(x)\gamma \theta(y))\\
&\geq\min\{\mu(\theta(x)),\mu(\theta(y))\}\\
&=\min\{\mu\lbrack\theta](x),\mu\lbrack\theta](y)\}.
\end{align*}
Thus $\mu\lbrack\theta]$ is a fuzzy subsemigroup of $S.$ Similarly we can prove the other case also.
\end{proof}

\begin{proposition}
If $\mu$ is a fuzzy subsemigroup$($fuzzy bi-ideal$)$ of a $\Gamma$-semigroup $S,$ then so is $\mu^{\alpha}$ for every real number $\alpha\geq 0,$ where $\mu^{\alpha}$ is defined by $\mu^{\alpha}(x)=(\mu(x))^{\alpha}$ for all $x\in S.$
\end{proposition}

\begin{proof}
Let $x,y\in S,\gamma\in\Gamma.$ Without any loss of generality, suppose $\mu(x)\geq\mu(y).$ Then $\mu^{\alpha}(x)\geq\mu^{\alpha}(y)$ and $\mu(x\gamma y)\geq\min\{\mu(x),\mu(y)\}=\mu (y).$ Then
\begin{align*}
\mu^{\alpha}(x\gamma y)&=(\mu(x\gamma y))^{\alpha}\geq(\mu(y))^{\alpha}=\mu^{\alpha}(y)\\
&=\min\{\mu^{\alpha}(x),\mu^{\alpha}(y)\}.
\end{align*}
Consequently, $\mu^{\alpha}$ is a fuzzy subsemigroup of $S.$ Similarly we can prove the other case also.
\end{proof}

\begin{theorem}
A non-empty fuzzy subset $\mu$ of a $\Gamma$-semigroup $S$ is a fuzzy subsemigroup of $S$ if and only if $\mu\circ\mu\subseteq\mu.$
\end{theorem}

\begin{proof}
Suppose $\mu\circ\mu\subseteq\mu.$ Then for $x,y\in S$ and $\gamma\in\Gamma$ we obtain $\mu(x\gamma y)\geq(\mu\circ\mu)(x\gamma y)\geq\min\{\mu(x),\mu(y)\}.$ So $\mu$ is a fuzzy subsemigroup of $S.$

Conversely, suppose $\mu$ is a fuzzy subsemigroup of $S$ and $x\in S.$  Suppose there exist $y,z\in S$ and $\gamma\in\Gamma$ such that $x=y\gamma z$ then $(\mu\circ\mu)(x)\neq 0.$ By hypothesis, $\mu(y\gamma z)\geq\min\{\mu(y),\mu(z)\}.$ Hence
\begin{align*}
\mu(y\gamma z)&\geq\underset{x=y\gamma z}{\sup}\min\{\mu(y),\mu(z)\}=(\mu\circ\mu)(x).
\end{align*}
Again if there does not exist $y,z\in S$ and $\gamma\in\Gamma$ such that $x=y\gamma z$ then $(\mu\circ\mu)(x)=0\leq\mu(x).$ Consequently, $\mu\circ\mu\subseteq\mu.$ This completes the proof.
\end{proof}

\begin{theorem}
In a $\Gamma$-semigroup $S$ for a non-empty fuzzy subset $\mu$ of $S$ the following are equivalent: $(1)$ $\mu$ is a fuzzy bi-ideal of $S,$ $(2)$ $\mu\circ\mu\subseteq\mu $ and $\mu\circ\chi\circ\mu\subseteq\mu ,$ where $\chi$ is the characteristic function of $S$.
\end{theorem}

\begin{proof}
Suppose $(1)$ holds, $i.e.,\mu$ is a fuzzy bi-ideal of $S.$ Then $\mu$ is a fuzzy subsemigroup of $S.$ So by Theorem $3.12$, $\mu\circ\mu\subseteq\mu.$ Let $a\in S.$ Suppose there exists $x,y,p,q\in S$ and $\beta,\gamma\in\Gamma$ such that $a=x\gamma y$ and $x=p\beta q.$ Since $\mu$ is a fuzzy bi-ideal of $S,$ we obtain $\mu(p\beta q\gamma y)\geq\min\{\mu(p),\mu(y)\}.$ Then
\begin{align*}
(\mu\circ\chi\circ\mu)(a)&=\underset{a=x\gamma y}{\sup}[\min\{(\mu\circ\chi)(x),\mu(y)\}]\\
&=\underset{a=x\gamma y}{\sup}[\min\{\underset{x=p\beta q}{\sup}\min\{\mu(p),\chi(q)\}\},\mu(y)]\\
&=\underset{a=x\gamma y}{\sup}[\min\{\underset{x=p\beta q}{\sup}\min\{\mu(p),1\}\},\mu(y)]\\
&=\underset{a=x\gamma y}{\sup}\{\min\{\mu(p),\mu(y)\}\}\\
&\leq\mu(p\beta q\gamma y)\\
&=\mu(x\gamma y)=\mu(a).
\end{align*}
So we have $(\mu\circ\chi\circ\mu)\subseteq\mu.$ Otherwise $(\mu\circ\chi\circ\mu)(a)=0\leq\mu(a).$ Thus
$(\mu\circ\chi\circ\mu)\subseteq\mu.$

Conversely, let us assume that $(2)$ holds. Since $\mu\circ\mu\subseteq\mu,$ so $\mu$ is a fuzzy subsemigroup of $S.$ Let $x,y,z\in S$ and $\beta,\gamma\in\Gamma$ and $a=x\beta y\gamma z.$ Since $\mu\circ\chi\circ\mu\subseteq\mu,$ we have
\begin{align*}
\mu(x\beta y\gamma z)&=\mu(a)\\
&\geq(\mu\circ\chi\circ\mu)(a)\\
&=\underset{a=x\beta y\gamma z}{\sup}[\min\{(\mu\circ\chi)(x\beta y),\mu(z)\}]\\
&\geq\min\{\mu\circ\chi)(p),\mu(z)\}(\text{let }p=x\beta y)\\
&=\min[\underset{p=x\beta y}{\sup}\min\{\mu(x),\chi(y)\},\mu(z)]\\
&\geq\min[\min\{\mu(x),1\},\mu(z)]\\
&=\min\{\mu(x),\mu(z)\}.
\end{align*}
Hence $\mu$ is a fuzzy bi-ideal of $S.$ This completes the proof.
\end{proof}

\section{Fuzzy Ideal and Fuzzy Bi-ideal}

\begin{definition}
A $\Gamma$-semigroup $S$ is called left$($right$)$ duo if every left$($resp. right$)$ ideal of $S$ is a two sided ideal of $S.$
\end{definition}

\begin{definition}
A $\Gamma$-semigroup $S$ is called duo if it is left and right duo.
\end{definition}

\begin{definition}
A $\Gamma$-semigroup $S$ is called fuzzy left$($right$)$ duo if every fuzzy left$($resp. right$)$ ideal of $S$ is a fuzzy ideal of $S.$
\end{definition}

\begin{definition}
A $\Gamma$-semigroup $S$ is called fuzzy duo if it is fuzzy left and fuzzy right duo.
\end{definition}

\begin{definition}
\cite{U}A $\Gamma$-semigroup $S$ is called left$($right$)$ regular if, for each $a\in S,$ there exist $x\in S$ and $\alpha,\beta\in\Gamma $ such that $a=x\alpha a\beta a($resp. $a=a\beta a\alpha x).$
\end{definition}

\begin{definition}
\cite{D1}A $\Gamma$-semigroup $S$ is called regular if, for each $a\in S,$ there exist $x\in S$ and $\alpha,\beta\in\Gamma$ such that $a=a\alpha x\beta a.$
\end{definition}

\begin{theorem}
In a regular left duo$($right duo, duo$)$ $\Gamma$-semigroup $S,$ following are equivalent: $(1)$ $\mu$ is a fuzzy right ideal$($resp. fuzzy left ideal, fuzzy ideal$)$ of $S,$ $(2) $ $\mu$ is a fuzzy bi-ideal of $S.$
\end{theorem}

\begin{proof}
Let $\mu$ be a fuzzy right ideal of $S$ and let $x,y,z\in S,\alpha,\beta\in\Gamma.$ Then
\begin{align*}
\mu(x\alpha y\beta z)=\mu(x\alpha (y\beta z))\geq\mu(x)\geq\min\{\mu(x),\mu(z)\}.
\end{align*}
Hence $\mu$ is a fuzzy bi-ideal of $S.$ Similarly we can prove the other cases.

Conversely, let $\mu$ be a fuzzy bi-ideal of $S$ and $x,y\in S,\gamma\in\Gamma .$ Then $x\gamma y\in S.$ Since $S$ is regular and left duo in view of the fact that $S\Gamma x$ is a left ideal we obtain, $x\gamma y\in(x\Gamma S\Gamma x)\Gamma S\subseteq x\Gamma S\Gamma x.$ This implies that there exist elements $z\in S,\alpha,\beta\in\Gamma$
such that $x\gamma y=x\alpha z\beta x.$ Then
\begin{align*}
\mu(x\gamma y)=\mu(x\alpha z\beta x)\geq\min\{\mu(x),\mu(x)\}=\mu(x).
\end{align*}
Hence $\mu$ is a fuzzy right ideal of $S.$ similarly we can prove the other cases also.
\end{proof}

\begin{theorem}
In a regular left duo$($right duo, duo$)$ $\Gamma$-semigroup $S,$ following are equivalent: $(1)$ $\mu$ is a fuzzy bi-ideal of $S,$ $(2)$ $\mu$ is a fuzzy $(1,2)$-ideal of $S.$
\end{theorem}

\begin{proof}
Let $\mu$ be a fuzzy bi-ideal of $S$ and $x,\omega ,y,z\in S,\alpha,\beta,\gamma\in\Gamma.$ Then
\begin{align*}
\mu(x\alpha\omega\beta(y\gamma z))&=\mu((x\alpha\omega\beta y)\gamma z)\\
&\geq\min\{\mu(x\alpha\omega\beta y),\mu(z)\}\\
&\geq\min\{\min\{\mu(x),\mu(y)\},\mu(z)\}\\
&=\min\{\mu(x),\mu(y),\mu(z)\}.
\end{align*}
Hence $\mu$ is a fuzzy $(1,2)$-ideal of $S.$

Conversely, let $S$ be a regular and left duo $\Gamma $-semigroup and $\mu$ be a fuzzy $(1,2)$-ideal of $S.$ Let $x,\omega,y\in S,\alpha,\delta\in\Gamma.$ Since $S$ is regular and left duo, we have $x\alpha\omega\in(x\Gamma S\Gamma x)\Gamma S\subseteq x\Gamma S\Gamma x,$ which implies that $x\alpha\omega=x\beta s\gamma x$ for some $s\in S$ and $\beta,\gamma\in\Gamma.$ Then
\begin{align*}
\mu(x\alpha\omega\delta y)&=\mu((x\beta s\gamma x)\delta y)\\
&=\mu(x\beta s\gamma(x\delta y))\\
&\geq\min\{\mu(x),\mu(x),\mu(y)\}\\
&=\min\{\mu(x),\mu(y)\}.
\end{align*}
Hence $\mu$ is a fuzzy bi-ideal of $S.$
\end{proof}

\begin{theorem}
For a regular $\Gamma$-semigroup $S$ the following conditions are equivalent: $(1)$ $S$ is left duo$($right duo, duo$),$ $(2)$ $S$ is fuzzy left duo$($fuzzy right duo, fuzzy duo$)$.
\end{theorem}

\begin{proof}
Let us assume that $S$ is left duo. Let $\mu$ be any fuzzy left ideal of $S$ and $a,b\in S,\gamma\in\Gamma.$ Then, since the left ideal $S\Gamma a$ is a two-sided ideal ideal of $S,$ and since $S$ is regular, we have $a\gamma b\in (a\Gamma S\Gamma a)\Gamma b\subseteq(S\Gamma a)\Gamma S\subseteq S\Gamma a.$ This implies that there exist elements $x\in S,\alpha\in\Gamma$ such that $a\gamma b=x\alpha a.$ Then, since $\mu$ is a fuzzy left ideal of $S,$ $\mu(x\alpha a)\geq\mu(a)$ and so $\mu(a\gamma b)\geq\mu(a).$ Hence $\mu$ is a fuzzy two-sided ideal of $S.$ Hence we deduce that $S$ is fuzzy left duo. Hence $(1)$ implies $(2).$ Similar for other cases.

Conversely, let us assume that $S$ is a fuzzy left duo. Let $A$ be any left ideal of $S.$ Then it follows from Theorem $3.1\cite{S2},$ that the characteristic function $\mu_{A}$ of $A$ is a fuzzy left ideal of $S.$ Thus by assumption it is a fuzzy ideal of $S.$ Since $A$ is non-empty, it follows from Theorem $3.1\cite{S2},$ that $A$ is an ideal
of $S.$ Therefore we obtain that $S$ is left duo, this completes the proof. Similarly we can prove the other cases.
\end{proof}

\begin{theorem}
For a regular $\Gamma$-semigroup $S$ the following conditions are equivalent: $(1)$ Every bi-ideal of $S$ is a right ideal$($left ideal, two-sided ideal$)$ of $S,$ $(2)$ every fuzzy bi-ideal of $S$ is a fuzzy right ideal$($ resp. fuzzy left ideal, fuzzy two-sided ideal$)$ of $S$.
\end{theorem}

\begin{proof}
Let us assume that every bi-ideal of $S$ is a right ideal. Let $\mu$ be a fuzzy bi-ideal of $S$ and $a,b\in S,\gamma\in\Gamma.$ Then $a\Gamma S\Gamma a$ is a bi-ideal of $S.$ Then by hypothesis $a\Gamma S\Gamma a$ is a right ideal of $S.$ Again since $S$ is regular, $a\gamma b\in(a\Gamma S\Gamma a)\Gamma S\subseteq a\Gamma S\Gamma a.$ This
implies that there exist elements $x\in S$ and $\alpha,\beta\in\Gamma$ such that $a\gamma b=a\alpha x\beta a.$ Then, since $\mu$ is a fuzzy bi-ideal of $S,$ we have $\mu(a\gamma b)=\mu(a\alpha x\beta a)\geq\min\{\mu(a),\mu(a)\}=\mu(a).$ Hence $\mu$ is a fuzzy right ideal of $S.$

Conversely, let us assume that every fuzzy bi-ideal is a fuzzy right ideal. Let $A$ be any bi-ideal of $S.$ Then it follows from Theorem $3.5,$ that the characteristic function $\mu_{A} $ of $A$ is a fuzzy bi-ideal of $S$ and consequently by hypothesis it is a fuzzy right ideal of $S.$ Then, since $A$ is non-empty, it follows from Theorem $3.1\cite{S2},$ that $A$ is a right ideal of $S.$ Similarly we can prove all other cases.
\end{proof}

\begin{definition}
\cite{U} A $\Gamma$-semigroup $S$ is called left-zero$($right zero$)$ if $x\gamma y=x(x\gamma y=y)$ $\forall x,y\in S,\forall\gamma\in\Gamma.$
\end{definition}

\begin{definition}
\cite{U} An element $e$ in a $\Gamma$-semigroup $S$ is called idempotent if $e\gamma e=e$ for some $\gamma\in\Gamma.$
\end{definition}

\begin{proposition}
For a left-zero$($right zero$)$ $\Gamma$-semigroup $S$ every fuzzy left$($resp. fuzzy right$)$ ideal is a constant function.
\end{proposition}

\begin{proof}
Let $S$ be a left zero $\Gamma$-semigroup and $\mu$ a fuzzy left ideal of $S.$ Let $x,y\in S.$ Then $x\gamma y=x$ and $y\gamma x=y \forall\gamma\in\Gamma.$ Then $\mu(x)=\mu(x\gamma y)\geq\mu(y).$ Again $\mu(y)=\mu(y\gamma x)\geq\mu(x).$ So $\mu(x)=\mu(y)$ $\forall x,y\in S.$ Hence every fuzzy left ideal is a constant function. Similarly we can prove the other case also.
\end{proof}

\begin{theorem}
For a regular $\Gamma$-semigroup $S$ the following conditions are equivalent: $(1)$ The set of all idempotent elements of $S$ forms a left zero$($resp. right zero$)$ subsemigroup of $S,$ $(2)$ for every fuzzy left$($resp. fuzzy right$)$ ideal $\mu$ of $S,\mu(e)=\mu(f)$ for all idempotent elements $e,f\in S.$
\end{theorem}

\begin{proof}
Let $E_{S}$ denote the set of all idempotent elements of $S$. Suppose $E_{S}$ is a left zero subsemigroup of $S.$ Let $e,f\in E_{S}$ and $\mu$ be a fuzzy left ideal of $S.$ Since, $S$ is left zero, then $e\gamma f=e$ and $f\gamma e=f$
$\forall\gamma\in\Gamma.$ Now, we have $\mu(e)=\mu(e\gamma f)\geq\mu(f)=\mu(f\gamma e)\geq\mu(e).$ Hence $\mu(e)=\mu(f).$

Conversely, let every fuzzy left ideal $\mu$ of $S$ satisfies the equality in $(2).$ Since $S$ is regular, for $a\in S$, there exists $x\in S$ and $\alpha,\beta\in\Gamma$ such that $a=a\alpha x\beta a\Rightarrow a\alpha x=(a\alpha x)\beta(a\alpha x).$ So $a\alpha x\in E_{S}$. So $E_{S}$ is non-empty. Let $e,f\in E_{S},\gamma\in\Gamma.$ Then from Theorem $3.1\cite{S2},$ the characteristic function $\mu_{L[f]}$ of the left ideal $L[f]$ of $S$ generated by $f,$ is a fuzzy left ideal of $S.$ Then by hypothesis $\mu_{L[f]}(e)=\mu_{L[f]}(f)=1$ and so $e\in L[f]=S\Gamma f.$ Then for some $x\in S,\alpha,\beta\in\Gamma$ we obtain $e=x\alpha f=x\alpha f\beta f=e\beta f.$ Hence $E_{S}$ is left zero. Now $e\gamma f=e($since $E_{S}$ is left zero$)=e\delta e($for some $\delta\in\Gamma)=(e\gamma f)\delta(e\gamma f).$ Consequently, $e\gamma f\in E_{S}.$ Hence $E_{S}$ is a left zero subsemigroup of $S.$ Similarly, we can prove the other case also.
\end{proof}

In view of the above theorem we obtain the following corollary.

\begin{corollary}
For an idempotent $\Gamma$-semigroup $S$ the following conditions are equivalent: $(1)$ $S$ is left zero$($resp. right zero$),$ $(2)$ for every fuzzy left$($resp. fuzzy right$)$ ideal $\mu$ of $S,\mu (e)=\mu (f)\forall e,f\in S.$
\end{corollary}

\begin{definition}
\cite{U}A $\Gamma$-semigroup $S$ is called intra-regular if, for each $a\in S,$ there exist $x,y\in S$ and $\alpha,\beta,\gamma\in\Gamma$ such that $a=x\alpha a\beta a\gamma y.$
\end{definition}

\begin{theorem}
For a $\Gamma$-semigroup $S$ the following conditions are equivalent: $(1)$ $S$ is intra-regular, $(2)$ for every fuzzy ideal $\mu$ of $S,\mu(a)=\mu(a\beta a)$ $\forall a\in S$ and for some $\beta\in\Gamma.$
\end{theorem}

\begin{proof}
Let us assume that $(1)$ holds. Let $\mu$ be a fuzzy ideal of $S$ and $a$ be any element of $S.$ Then there exist $x,y\in S$ and $\alpha,\beta,\gamma\in\Gamma $ such that $a=x\alpha a\beta a\gamma y.$ Since $\mu$ is a fuzzy ideal of $S,$ we have $\mu(a)=\mu(x\alpha a\beta a\gamma y)\geq\mu(x\alpha a\beta a)\geq\mu(a\beta a)\geq\min\{\mu(a),\mu(a)\}=\mu (a).$ So $\mu(a)=\mu(a\beta a).$

Conversely, let us assume that $(2)$ holds. Then it follows from Theorem $3.1\cite{S2},$ that the characteristic function $\mu_{J[a\beta a]}$ of the ideal $J[a\beta a]$ of $S$ generated by $a\beta a,$ is a fuzzy ideal of $S.$ Since $a\beta a\in J[a\beta a],$ we have $\mu_{J[a\beta a]}(a)=\mu_{J[a\beta a]}(a\beta a)=1.$ This implies that $a\in J[a\beta a].$ This proves that $S$ is intra-regular.
\end{proof}

\begin{theorem}
For a $\Gamma$-semigroup $S$ the following conditions are equivalent: $(1)$ $S$ is left regular$($resp. right regular$),$ $(2)$ for every fuzzy left ideal$($resp. fuzzy right ideal$)$ $\mu$ of $S,\mu(a)=\mu(a\gamma a)$ $\forall a\in S$ and for some $\gamma\in\Gamma.$
\end{theorem}

\begin{proof}
Let us assume that $S$ is left regular. Let $\mu$ be any fuzzy left ideal of $S$ and $a\in S.$ Then there exist $x\in S$ and $\alpha,\gamma\in\Gamma $ such that $a=x\alpha a\gamma a$ and so $\mu(a)=\mu(x\alpha a\gamma a)\geq\mu(a\gamma a)\geq\min\{\mu(a),\mu(a)\}=\mu(a).$ Consequently $\mu(a)=\mu(a\gamma a).$

Conversely, let us assume that for every fuzzy left ideal $\mu$ of $S$ the equality in $(2)$ holds. Let $a$ be any element of $S.$ Then it follows from Theorem $3.1\cite{S2},$ that characteristic function $\mu_{L[a\gamma a]}$ of the left ideal $L[a\gamma a]$ of $S$ generated by $a\gamma a,$ is a fuzzy left ideal of $S.$ Since $a\gamma a\in L[a\gamma a],$ we have $\mu_{L[a\gamma a]}(a)=\mu_{L[a\gamma a]}(a\gamma a)=1.$ This implies that $a\in L[a\gamma a]=\{a\gamma a\}\cup S\Gamma a\Gamma a.$ This proves that $S$ is left regular.  Similarly we can prove the result for fuzzy right ideals.
\end{proof}

\begin{proposition}
Suppose $S$ is both regular and intra regular $\Gamma$-semigroup. Then $(1)$ $\mu_{1}\circ\mu_{2}\supset\mu_{1}\cap\mu_{2}.$ $(2)$ $(\mu_{1}\circ\mu_{2})\cap(\mu_{2}\circ\mu_{1})\supset\mu_{1}\cap \mu_{2}$ where $\mu_{1},\mu_{2}$ are fuzzy bi-ideals of $S.$
\end{proposition}

\begin{proof}
Let $a\in S.$ Then there exist $x,y,z\in S$ and $\gamma_{1},\gamma_{2}\gamma_{3},\gamma_{4},\gamma_{5}\in\Gamma $ such that $a=a\gamma_{1}x\gamma_{2}a=a\gamma_{1}x\gamma_{2}a\gamma_{1}x\gamma_{2}a=a\gamma_{1}x\gamma_{2}(y\gamma_{3}a\gamma _{4}a\gamma_{5}z)\gamma_{1}x\gamma_{2}a=(a\gamma_{1}x\gamma_{2}y\gamma_{3}a)\gamma_{4}(a\gamma_{5}z\gamma_{1}x\gamma _{2}a).$ Since $\mu_{1},\mu_{2}$ are both fuzzy bi-ideals of $S,$ we deduce that $\mu_{1}(a\gamma_{1}x\gamma_{2}y\gamma _{3}a)\geq\mu_{1}(a)$ and $\mu_{2}(a\gamma_{5}z\gamma_{1}x\gamma_{2}a)\geq\mu_{2}(a).$ Then
\begin{align*}
(\mu_{1}\circ\mu_{2})(a)&=\underset{a=p\gamma_{2}q}{\sup}[\min\{\mu_{1}(p),\mu_{2}(q)\}:p,q\in S;\gamma_{2}\in\Gamma]\\
&\geq\min\{\mu_{1}(a\gamma_{1}x\gamma_{2}y\gamma_{3}a),\mu_{2}(a\gamma_{5}z\gamma_{1}x\gamma_{2}a)\}
\end{align*}
\begin{align*}
\geq\min\{\mu_{1}(a),\mu_{2}(a)\}&=(\mu_{1}\cap\mu_{2})(a).
\end{align*}
Hence $\mu_{1}\circ\mu_{2}\supset\mu_{1}\cap\mu_{2}$. Similarly we can show that $\mu_{2}\circ\mu_{1}\supset\mu_{1}\cap\mu_{2}.$ Therefore $(\mu_{1}\circ\mu_{2})\cap(\mu_{2}\circ\mu_{1})\supset\mu_{1}\cap\mu_{2}$. This completes the proof.
\end{proof}

\begin{definition}
A $\Gamma$-semigroup $S$ is said to be left$($right$)$ simple if $S$ has no proper left$($resp. right$)$ ideals.
\end{definition}

\begin{definition}
If a $\Gamma$-semigroup $S$ has no proper ideals, then we say that $S$ is simple.
\end{definition}

\begin{definition}
A $\Gamma$-semigroup $S$ is said to be fuzzy left$($fuzzy right$)$ simple if every fuzzy left$($resp. fuzzy right$)$ ideal of $S$ is a constant function.
\end{definition}

\begin{definition}
A $\Gamma$-semigroup $S$ is said to be fuzzy simple if every fuzzy ideal of $S$ is a constant function.
\end{definition}

\begin{theorem}
For a $\Gamma$-semigroup $S$ the following conditions are equivalent: $(1)$ $S$ is left simple$($resp. right simple,
simple$),$ $(2)$ $S$ is fuzzy left simple$($resp. fuzzy right simple, fuzzy simple$)$.
\end{theorem}

\begin{proof}
Let us assume that $S$ is left simple. Let $\mu$ be any fuzzy left ideal of $S$ and $a,b\in S.$ Then there exist $x,y\in S$ and $\alpha,\beta\in\Gamma$ such that $b=x\alpha a$ and $a=y\beta b$ and so we obtain $\mu(a)=\mu(y\beta b)\geq\mu(b)=\mu(x\alpha a)\geq\mu(a).$ Consequently, $\mu (a)=\mu(b)$ and so $\mu$ is a constant function. Hence $S$ is fuzzy left simple.

Conversely, let us assume that $S$ is fuzzy left simple and let $A$ be any left ideal of $S.$ Then by Theorem $3.1\cite{S2},$ $\mu_{A}$ is a fuzzy left ideal of $S$ and hence a constant function. Since $A$ is non-empty, the constant is $1.$ So every element of $S$ is in $A$ and so $S$ is left simple. Similarly we can prove the other cases.
\end{proof}

\begin{theorem}
Let $S$ be a left$($right$)$ simple $\Gamma$-semigroup, then every fuzzy bi-ideal of $S$ is a fuzzy right ideal$($resp. fuzzy left ideal$)$ of $S.$
\end{theorem}

\begin{proof}
Let $S$ be a left simple $\Gamma$-semigroup. Let $\mu$ be any fuzzy bi-ideal of $S$ and $a,b\in S.$ Then there exist $x\in S,\gamma\in\Gamma $ such that $b=x\gamma a$ and $\mu(a\alpha b)=\mu(a\alpha x\gamma a)\geq\min\{\mu(a),\mu(a)\}=\mu(a)\forall\alpha\in\Gamma.$ Hence $\mu$ is a fuzzy right ideal of $S.$ Similarly we can prove the other case also.
\end{proof}

\section{Fuzzy Quasi Ideal}

\begin{definition}
A non-empty fuzzy subset $\mu$ of a $\Gamma$-semigroup $S$ is called a fuzzy quasi ideal of $S$ if $(\mu\circ\chi)\cap(\chi\circ\mu)\subseteq\mu,$ where $\chi$ is the characteristic function of $S.$
\end{definition}

\begin{proposition}
Any fuzzy one sided ideal of a $\Gamma$-semigroup $S$ is a fuzzy quasi ideal of $S$ and any fuzzy quasi ideal of $S$ is a fuzzy bi-ideal of $S.$
\end{proposition}

\begin{proof}
Let $\mu$ be any fuzzy left ideal of $S.$ Then $\chi\circ\mu\subseteq\mu$. Therefore $(\mu\circ\chi)\cap(\chi\circ\mu)\subseteq\chi\circ\mu\subseteq\mu$. Therefore $\mu$ is fuzzy quasi ideal of $S.$
Again let $\mu$ be a fuzzy quasi ideal of $S.$ Then $(\mu\circ\chi)\cap(\chi\circ\mu)\subseteq\mu$. Now $\mu\circ\chi\circ\mu\subseteq\mu\circ\chi$ and $\mu\circ\chi\circ\mu\subseteq\chi\circ\mu$. Therefore $\mu\circ\chi\circ\mu\subseteq(\mu\circ\chi)\cap(\chi\circ\mu)\subseteq\mu$. Consequently, $\mu$ is fuzzy bi ideal of $S.$
\end{proof}

\begin{proposition}
Any fuzzy quasi ideal of a $\Gamma$-semigroup $S$ can be expressed as intersection of a fuzzy right ideal and a fuzzy left ideal of $S$ and conversely intersection of a fuzzy right ideal and a fuzzy left ideal is fuzzy quasi ideal .
\end{proposition}

\begin{proof}
Let $\mu$ be a fuzzy quasi ideal of a $\Gamma$-semigroup $S.$ Then $(\chi\circ\mu)\cup\mu$ is a fuzzy left ideal and $\mu\cup(\mu\circ\chi)$ is a fuzzy right ideal of $S.$ Clearly $\mu\subseteq\mu\cup(\mu\circ\chi)$ and $\mu\subseteq (\chi\circ\mu)\cup\mu$. Therefore
\begin{align*}
\mu &\subseteq(\mu\cup(\mu\circ\chi))\cap((\chi\circ\mu)\cup\mu)\\ &=(\mu\cap(\chi\circ\mu))\cup(\mu\cap\mu)\cup((\mu\circ\chi)\cap(\chi\circ\mu))\cup((\mu\circ\chi)\cap\mu)\\ &\subseteq\mu\cup\mu\cup\mu\cup\mu =\mu.
\end{align*}
Therefore $\mu=(\mu\cup(\mu\circ\chi))\cap((\chi\circ\mu)\cup\mu)$. Let $\mu$ and $\sigma$ be a fuzzy right ideal and a fuzzy left ideal of the $\Gamma$-semigroup $S$ respectively. Let $\chi$ be the characteristic function of $S.$ Since $\mu$ and $\sigma$ are fuzzy right ideal and fuzzy left ideal of $S,$ $\mu\circ\chi\subseteq\mu$ and $\chi\circ\sigma\subseteq\sigma(cf.$ Theorem $4.2\cite{S2}).$

Now, $((\mu\cap\sigma)\circ\chi)\cap(\chi\circ(\mu\cap\sigma))\subseteq(\mu\circ\chi)\cap(\chi\circ\sigma)\subseteq\mu\cap \sigma.$ Thus $\mu\cap\sigma$ is a fuzzy quasi ideal of $S.$ This completes the proof.
\end{proof}

In view of Proposition $5.3$ and Theorem $4.7\cite{S2}$ we have the following corollary.

\begin{corollary}
Let $\mu$ and $\sigma$ be a fuzzy right ideal and a fuzzy left ideal of a regular $\Gamma$-semigroup $S,$ respectively. Then $\mu\circ\sigma $ is a fuzzy quasi ideal of $S.$
\end{corollary}

In view of Proposition $5.2$ and Corollary $5.4$ we obtain the following proposition.

\begin{proposition}
Let $\mu$ be a fuzzy right ideal and $\sigma$ be a fuzzy left ideal of a regular $\Gamma$-semigroup $S.$ Then $\mu\circ\sigma $ is a fuzzy bi-ideal of $S.$
\end{proposition}

\begin{theorem}
Let $Q$ be a non-empty subset of a $\Gamma$-semigroup $S.$ Then $Q$ is quasi ideal of $S$ if and only if $\mu_{Q}($the characteristic function of $Q)$ is a fuzzy quasi ideal of $S.$
\end{theorem}

\begin{proof}
Suppose $Q$ is a quasi ideal of $S$ and $\chi$ be the characteristic function of $S.$ Let $a$ be any element of $S.$ If $a\in Q,$ then $((\mu_{Q}\circ\chi )\cap(\chi\circ\mu_{Q}))(a)\leq 1=\mu_{Q}(a). $

If $a\notin Q,$ then since $Q$ is a quasi ideal of $S,i.e.,Q\Gamma S\cap S\Gamma Q\subseteq Q$, $a\notin Q\Gamma S\cap
S\Gamma Q\subseteq Q.$ Then three cases may arise:
\begin{align*}
&\text{Case-(1) }a\notin Q\Gamma S, a\in S\Gamma Q,\\
&\text{Case-(2) }a\in Q\Gamma S, a\notin S\Gamma Q,\\
&\text{Case-(3) }a\notin Q\Gamma S, a\notin S\Gamma Q.
\end{align*}

Case-$(1):$ Let $a\notin Q\Gamma S, a\in S\Gamma Q.$ If $a=x\gamma y,$ then $x\notin Q.$ Then
\begin{align*}
((\mu_{Q}\circ\chi)\cap(\chi\circ\mu_{Q}))(a)&=\min[\underset{a=x\gamma y}{\sup }\{\min \{\mu _{Q}(x),\chi(y)\}\},\underset{a=u\theta v}{\sup }\{\min\{\chi(u),\mu_{Q}(v)\}\}]\\
&=0=\mu_{Q}(a).
\end{align*}
Therefore $(\mu \circ \chi)\cap (\chi \circ \mu )\subseteq \mu.$

Case-$(2):$ It is similar as Case-$(3).$

Case-$(3):$ $a\notin Q\Gamma S, a\notin S\Gamma Q.$ If $a=x\gamma y,$ then $x\notin Q$ and if $a=u\delta v,$ then $v\notin Q.$ Now
\begin{align*}
((\mu_{Q}\circ\chi)\cap(\chi\circ\mu_{Q}))(a)&=\min[\underset{a=x\gamma y}{\sup}\{\min\{\mu_{Q}(x),\chi(y)\}\},\underset{a=u\delta v}{\sup}\{\min\{\chi(u),\mu_{Q}(v)\}\}]\\
&=0=\mu_{Q}(a).
\end{align*}
So $(\mu\circ\chi)\cap(\chi\circ\mu)\subseteq\mu.$ Hence $\mu_{Q}$ is a fuzzy quasi ideal of $S.$

Conversely, let us suppose that $\mu_{Q}$ is a fuzzy quasi ideal of $S.$ Let $a\in Q\Gamma S\cap S\Gamma Q.$ Then there exist elements $s,t\in S,b,c\in Q$ and $\alpha,\beta\in\Gamma$ such that $a=b\alpha s=t\beta c.$ Then
\begin{align*}
(\mu_{Q}\circ \chi)(a)&=\underset{a=x\beta y}{\sup}[\min\{\mu_{Q}(x),\chi(y)\}]\\
&\geq\min\{\mu_{Q}(b),\chi(s)\}\\
&=\min\{1,1\}=1.
\end{align*}
Similarly, we have $(\chi\circ\mu_{Q})(a)=1.$ Since $(\mu_{Q}\circ\chi)\cap(\chi\circ\mu_{Q})\subseteq\mu_{Q}.$ So
\begin{align*}
\mu_{Q}(a)&\geq((\mu_{Q}\circ\chi )\cap(\chi\circ\mu_{Q}))(a)\\
&=\min\{(\mu_{Q}\circ\chi)(a),(\chi\circ\mu_{Q})(a)\}\\
&=\min \{1,1\}=1.
\end{align*}
Thus $a\in Q$ and so $Q\Gamma S\cap S\Gamma Q\subseteq Q.$ Therefore, $Q$ is a quasi ideal of $S.$
\end{proof}

\begin{theorem}
Let $\mu$ be a non-empty fuzzy subset of a $\Gamma$-semigroup $S.$ Then $\mu$ is a fuzzy quasi ideal of $S$ if and only if $\mu_{t}=\{x\in S:\mu(x)\geq t\}$ is a quasi ideal of $S,$ where $t\in Im(\mu).$
\end{theorem}

\begin{proof}
Let $\mu$ be a fuzzy quasi ideal of $S.$ Let $t\in Im(\mu)$,then $\mu_{t}$ is nonempty. Let $a\in S\Gamma \mu _{t}\cap \mu_{t}\Gamma S.$ Then there exist elements $r\in S,b\in\mu_{t}$ and $\alpha\in\Gamma $ such that $a=b\alpha r.$ Then
\begin{align*}
(\mu\circ \chi )(a)=\underset{a=x\beta y}{\sup }[\min \{\mu (x),\chi(y)\}]\geq \min \{\mu(b),\chi(r)\}\geq\min\{t,1\}=t.
\end{align*}
Similarly, we have $(\chi \circ \mu )(a)\geq t.$ Then
\begin{align*}
((\chi \circ \mu )\cap(\mu \circ \chi ))(a)=\min \{(\chi \circ \mu )(a),(\mu \circ \chi)(a)\}\geq t..
\end{align*}
By hypothesis $(\chi\circ\mu)\cap(\mu\circ\chi)\subseteq\mu .$ Hence $\mu(a)\geq ((\chi \circ \mu )\cap (\mu \circ \chi ))(a)\geq t.$ Thus $a\in \mu _{t}.$ Hence $S\Gamma \mu _{t}\cap \mu _{t}\Gamma S\subseteq \mu _{t}.$ Consequently $\mu _{t}$ is a quasi ideal of $S.$

Conversely, let $\mu_{t}$ is a quasi ideal of $S$ $\forall t\in Im(\mu).$ Then $\mu_{t}$ is nonempty and so $\mu$ is nonempty. If possible, let $((\mu \circ \chi )\cap (\chi \circ \mu))\not\subseteq \mu .$ Then there exists $p\in S$ such that $\mu(p)<((\mu\circ\chi)\cap(\chi\circ\mu))(p).$ Let $t_{2}$ be a
real number such that
\begin{align*}
\mu(p)<t_{2}<((\mu\circ\chi)\cap(\chi\circ\mu))(p)........\text{(1)}.
\end{align*}
If there not exist $x,y\in S,\gamma\in\Gamma$ such that $p=x\gamma y$ then $((\mu\circ\chi)\cap(\chi\circ\mu))(p)=0$ which means $\mu(p)<0,$ which is not possible. So, there exist $x,y\in S,\gamma\in\Gamma$ such that $p=x\gamma y.$ Now suppose for any $x,v\in S,$ with $p=x\gamma y, p=u\delta v$ for some $y,u\in S$ for some $\gamma,\delta\in\Gamma,$ $x\notin\mu_{t_{2}}$ and $v\notin\mu_{t_{2}}.$ Then $\mu(x)<t_{2}$ and $\mu(v)<t_{2}$ for all $x,v\in S$ with $p=x\gamma y$ and $p=u\delta v,$ $x,u\in S,\gamma,\delta\in\Gamma.$ This implies that $\underset{p=x\gamma y}{\sup}\mu(x)\leq t_{2}$ and $\underset{p=u\delta v}{\sup}\mu(v)\leq t_{2}.$ Therefore $\min\{\underset{p=x\gamma y}{\sup}\mu(x),\underset{p=u\delta v}{\sup }\mu(v)\}\leq t_{2},i.e.,((\mu\circ\chi)\cap(\chi\circ\mu))(p)\leq t_{2},$ this contradicts $(1).$ So there exist $x,v\in S$ with $p=x\gamma y,p=u\delta v$ for some $\gamma,\delta\in\Gamma,$ for some $y,u\in S$ such that $x\in\mu_{t_{2}}$ and $v\in\mu_{t_{2}}.$ Hence $x\gamma y\in\mu_{t_{2}}\Gamma S$ and $u\delta v\in S\Gamma\mu_{t_{2}}$ whence $p\in\mu_{t_{2}}\Gamma S\cap S\Gamma\mu_{t_{2}}.$ But by $(1), p\not\in\mu_{t_{2}}.$ This contradicts that $\mu_{t}$ is a quasi ideal of $S,$ for all $t\in Im(\mu)$. Hence $((\mu\circ\chi)\cap(\chi\circ\mu))\subseteq\mu.$ Consequently, $\mu$ is a fuzzy quasi ideal of $S.$
\end{proof}

\begin{proposition}
Let $f:R\rightarrow S$ be a homomorphism of $\Gamma$-semigroups. If $\lambda$ is a fuzzy left$($resp. right$)$ ideal of $S$, then $f^{-1}(\chi)\circ f^{-1}(\lambda)\subseteq f^{-1}(\lambda)$ $(resp.f^{-1}(\lambda)\circ f^{-1}(\chi)\subseteq f^{-1}(\lambda)),$ $\chi$ is the characteristic function of $S$ $($where $f^{-1}(\lambda)(r\gamma s):=\lambda(f(r\gamma s))$ for all $r,s\in R$ and $\gamma\in\Gamma),$ provided $f^{-1}(\lambda)$ is non-empty.
\end{proposition}

\begin{proof}
$f^{-1}(\lambda)$ is a fuzzy left ideal of $R(cf.$ Theorem $3.6\cite{S2}).$ Let $x\in R.$ Then if there exists $a,b\in R$ and $\alpha\in\Gamma$ such that $x=a\alpha b$ we obtain
\begin{align*}
(f^{-1}(\chi)\circ f^{-1}(\lambda))(x)&=\underset{x=u\gamma v}{\sup}[\min\{f^{-1}(\chi)(u),f^{-1}(\lambda)(v)\}]\\
&=\underset{x=u\gamma v}{\sup}[\min\{\chi(f(u)),\lambda(f(v))\}]\\
&=\underset{x=u\gamma v}{\sup}[\min\{1,f^{-1}(\lambda)(v)\}]\\
&\leq f^{-1}(\lambda)(u\gamma v)(\text{for all }u,v\in R,\text{ for all }\gamma\in\Gamma \text{ with }x=u\gamma v)\\
&=f^{-1}(\lambda )(x).
\end{align*}
Hence $f^{-1}(\chi)\circ f^{-1}(\lambda)\subseteq f^{-1}(\lambda).$ Similarly we can prove the other case also.
\end{proof}

In view of the above proposition we obtain the following proposition.

\begin{proposition}
Let $f:R\rightarrow S$ be a homomorphism of $\Gamma$-semigroups. If $\lambda$ is a fuzzy quasi ideal of $S$, then $f^{-1}(\lambda)$ is a fuzzy quasi ideal $R$, provided $f^{-1}(\lambda)$ is non-empty.
\end{proposition}

\begin{proof}
Let $\lambda$ be a fuzzy quasi ideal of $S$ and $\lambda=\mu\cap\nu,$ where $\mu$ is a fuzzy right ideal and $\nu$ is a fuzzy left ideal of $S(cf.$ Proposition $5.3).$ Let $a\in S.$ Then
\begin{align*}
f^{-1}(\lambda)(a)&=\lambda(f(a))\\
&=(\mu\cap\nu)(f(a))\\
&=\min\{\mu(f(a)),\nu(f(a)\}\\
&=\min\{f^{-1}(\mu)(a),f^{-1}(\nu)(a)\}\\
&=(f^{-1}(\mu)\cap f^{-1}(\nu))(a).
\end{align*}
Consequently, $f^{-1}(\lambda)=f^{-1}(\mu)\cap f^{-1}(\nu).$ By Proposition $3.6\cite{S2},$ $f^{-1}(\mu)$ is a fuzzy right ideal and $f^{-1}(\nu)$ is a fuzzy left ideal of $R.$ Hence $f^{-1}(\lambda)$ is a fuzzy quasi ideal of $R.$
\end{proof}

\begin{proposition}
Let $f:S\rightarrow R$ be a surjective homomorphism of $\Gamma$-semigroups. If $\lambda$ is a fuzzy left$($resp. right$)$ ideal of $S$, then $f(\chi)\circ f(\lambda)\subseteq f(\lambda)($resp. $f(\lambda)\circ f(\chi)\subseteq f(\lambda)),$ $\chi$ is the characteristic function of $S$ $($where $(f(\lambda ))(r)$ $\ :=\underset{f(s)=r}{\sup }\lambda (s)$ for $s\in S).$
\end{proposition}

\begin{proof}
Suppose $\lambda$ is a fuzzy left ideal of $S.$ Then $f(\lambda)$ is a fuzzy left ideal of $R\cite{S2}.$ For $r\in R$, $(f(\lambda))(r)=\underset{f(s)=r}{\sup}\lambda (s)$, so $f(\lambda)$ is non-empty. Then for $p\in R$,
\begin{align*}
(f(\chi)\circ f(\lambda))(p)&=\underset{p=u\gamma v}{\sup}[\min\{f(\chi)(u),f(\lambda)(v)\}]\\
&=\underset{p=u\gamma v}{\sup}[\min\{\underset{f(u^{^{\prime}})=u}{\sup}\chi(u^{^{\prime}}),f(\lambda)(v)\}]\\
&=\underset{p=u\gamma v}{\sup}[\min\{1,f(\lambda)(v)\}]\\
&\leq\underset{p=u\gamma v}{\sup}[\min\{1,f(\lambda)(u\gamma v)\}]\\
&=\underset{p=u\gamma v}{\sup}[\min\{1,f(\lambda)(p)\}]\\
&=f(\lambda)(p).
\end{align*}
Hence $f(\chi)\circ f(\lambda)\subseteq f(\lambda).$ Similarly we can prove the other case also.
\end{proof}

\begin{proposition}
Let $f:R\rightarrow S$ be a surjective homomorphism of $\Gamma$-semigroups. If $\lambda$ is a fuzzy quasi ideal of $R$, then $f(\lambda)$ is a fuzzy quasi ideal $S$.
\end{proposition}

\begin{proof}
Let $\lambda$ be a fuzzy quasi ideal of $R$. Then by Proposition $5.3,$ $\lambda=\mu\cap\nu,$ for some fuzzy right ideal $\mu$ and fuzzy left ideal $\nu$ of $R.$ Hence by Proposition $3.6\cite{S2},$ $f(\mu)$ is a fuzzy right ideal and $f(\nu)$ is a fuzzy left ideal of $S.$ Let $s\in S.$ Then
\begin{align*}
f(\lambda)(s)&=\underset{f(r)=s}{\sup}\lambda(r)(\text{for some }r\in R)\\
&=\underset{f(r)=s}{\sup}(\mu\cap\nu)(r)\\
&=\underset{f(r)=s}{\sup}\min\{\mu(r),\nu(r)\}\\
&=\min\{\underset{f(r)=s}{\sup}\mu(r),\underset{f(r)=s}{\sup}\nu(r)\}\\
&=\min\{f(\mu)(s),f(\nu)(s)\}\\
&=(f(\mu)\cap f(\nu))(s).
\end{align*}
Consequently, $f(\lambda)=f(\mu)\cap f(\nu).$ So, $f(\lambda)$ is the intersection of a fuzzy right ideal $f(\mu)$ and a fuzzy left ideal $f(\nu)$ of $S.$ Hence by Proposition $5.3,$ $f(\lambda) $ is a fuzzy quasi ideal $S$.
\end{proof}

\begin{lemma}
Let $S$ be a $\Gamma$-semigroup and $A,B\subseteq S.$ Then $(1)$ $A\subseteq B$ if and only if $\mu_{A}\subseteq\mu_{B}.$ $(2)$ $\mu_{A}\cap\mu_{B}=\mu_{A\cap B}.$ $(3)$ $\mu_{A}\circ\mu_{B}=\mu_{A\Gamma B},$ where $\mu_{A},\mu_{B}$ denote the characteristic functions of $A$ and $B$ respectively.
\end{lemma}

\begin{proof}
$(1)$ The proof follows by routine verification.\\

$(2)$ Let $a\in S.$ Suppose that $a\in A\cap B.$ Then $a\in A$ and $a\in B,$ which implies $\mu_{A}(a)=\mu_{B}(a)=1.$ Then
\begin{align*}
(\mu_{A}\cap\mu_{B})(a)&=\min\{\mu_{A}(a),\mu_{B}(a)\}\\
&=1=\mu_{A\cap B}(a).
\end{align*}
Again if, $a\notin A\cap B.$ Then $a\notin A$ or $a\notin B,$ which implies $\mu_{A}(a)=0$ or $\mu_{B}(a)=0.$ Then
\begin{align*}
(\mu_{A}\cap\mu_{B})(a)&=\min\{\mu_{A}(a),\mu_{B}(a)\}\\
&=0=\mu_{A\cap B}(a).
\end{align*}
Consequently $\mu_{A}\cap\mu_{B}=\mu_{A\cap B}.$\\

$(3)$ Let $a\in S.$ Suppose that $a\in A\Gamma B.$ Then $a=x\gamma y$ for some $x\in A,y\in B$ and $\gamma\in\Gamma.$ Then
\begin{align*}
(\mu_{A}\circ\mu_{B})(a)&=\underset{a=u\delta v}{\sup}\min\{\mu_{A}(u),\mu_{B}(v)\}\\
&\geq\min\{\mu_{A}(x),\mu_{B}(y)\}\\
&=\min\{1,1\}=1.
\end{align*}
So $(\mu_{A}\circ\mu_{B})(a)=1.$ Since $a\in A\Gamma B,\mu_{A\Gamma B}(a)=1.$ In the case, when $a\notin A\Gamma B$  then we have, $(\mu_{A}\circ\mu_{B})(a)=0=\mu_{A\Gamma B}(a).$ Thus we obtain $\mu_{A}\circ\mu_{B}=\mu_{A\Gamma B}.$

\end{proof}

The following theorem is the characterization of regular $\Gamma$-semigroup in terms of fuzzy quasi ideals.

\begin{theorem}
In a $\Gamma$-semigroup $S$ the following are equivalent: $(1)$ $S$ is regular, $(2)$ for every fuzzy right ideal $\mu$ and every fuzzy left ideal $\lambda$ of $S$ we have $\mu\circ\lambda=\mu\cap\lambda,$ $(3)$ for every fuzzy right ideal $\mu$ and every fuzzy left ideal $\lambda$ of $S$ we have $(a)$ $\mu\circ\mu=\mu,$ $(b)$ $\lambda\circ\lambda=\lambda,$ $(c)$ $\mu\circ\lambda$ is the fuzzy quasi ideal of $S,$ $(4)$ every fuzzy quasi ideal $\delta$ has the form $\delta=\delta\circ\chi_{S}\circ\delta,$ where $\chi_{S}$ is the characteristic function of $S.$
\end{theorem}

\begin{proof}
$(1)\Rightarrow (2):$ by Theorem $4.7(\cite{S2}).$\\

$(2)\Rightarrow (3):$ $(a)$ and $(b)$ can be verified by routine calculation.\\

$(c)$ From Proposition $5.3,$ it is clear that the intersection of any fuzzy right ideal and fuzzy left ideal of $S$ is a fuzzy quasi-ideal of $S.$ Hence by $(2)$ the result follows.\\

$(3)\Rightarrow (4):$ Let $(3)$ holds. Let $\delta$ be a fuzzy quasi ideal of $S.$ Then $\delta\circ\chi_{S}$ and $\chi_{S}\circ\delta$ are fuzzy right and left ideals of $S$ respectively.
Then by conditions $(b),(c)$ we obtain
\begin{align*}
(\delta\circ\chi_{S})\circ(\chi_{S}\circ\delta)=\delta\circ(\chi_{S}\circ\chi_{S})\circ\delta=\delta\circ\chi_{S}\circ\delta
\end{align*}
is a fuzzy quasi ideal of $S.$ Also we have
\begin{align*}
\delta\circ\chi_{S}\circ\delta\subseteq(\delta\circ\chi_{S})\cap(\chi_{S}\circ\delta)\subseteq\delta..............(A).
\end{align*}
Now
\begin{align*}
\delta&\subseteq\delta\cup(\delta\circ\chi_{S})\\
&=(\delta\cup(\delta\circ\chi_{S}))\circ(\delta\cup(\delta\circ\chi_{S}))(\text{by condition }(3)(a))\\
&=((\delta\cup(\delta\circ\chi_{S}))\circ\delta)\cup((\delta\cup(\delta\circ\chi_{S}))\circ(\delta\circ\chi_{S}))\\
&=((\delta\circ\delta)\cup(\delta\circ\chi_{S}\circ\delta))\cup((\delta\circ\delta\circ\chi_{S})
\cup((\delta\circ\chi_{S})\circ(\delta\circ\chi_{S})))\\
&\subseteq((\delta\circ\chi_{S})\cup(\delta\circ\chi_{S}))
\cup((\delta\circ\chi_{S})\cup(\delta\circ\chi_{S}))\\
&=\delta\circ\chi_{S}(\text{since }\delta\circ\delta
\subseteq\delta,\delta\circ\delta\subseteq\delta\circ\chi_{S}\text{ and }\delta\circ\chi_{S}\text{ is a fuzzy right ideal of }S).
\end{align*}
Similarly we obtain $\delta\subseteq\chi_{S}\circ\delta$. Then $\delta\subseteq(\delta\circ\chi_{S})\cap(\chi_{S}\circ\delta)$. Thus
\begin{align*}
\delta&=(\delta\circ\chi_{S})\cap(\chi_{S}\circ\delta)\\
&=((\delta\circ\chi_{S})\circ(\delta\circ\chi_{S}))\cap((\chi_{S}\circ\delta)\circ(\chi_{S}\circ\delta))(\text{by condition }(a),(b))\\
&=((\delta\circ\chi_{S}\circ\delta)\circ\chi_{S})\cap(\chi_{S}\circ(\delta\circ\chi_{S}\circ\delta))\\
&\subseteq\delta\circ\chi_{S}\circ\delta(\text{since }\delta\circ\chi_{S}\circ\delta\text{ is fuzzy quasi ideal of } S).........(B).
\end{align*}
By $(A)$ and $(B)$ we obtain $\delta=\delta\circ\chi_{S}\circ\delta$.

$(4)\Rightarrow (1):$ Let $(4)$ holds. Let $x\in S.$ Let us consider the quasi ideal $Q[x]=\{x\}\cup(x\Gamma S\cap S\Gamma x)$ generated by $x.$ Then the characteristic function $\chi_{Q[x]}$ of $Q[x]$ is a fuzzy quasi ideal of $S.$ Then
\begin{align*}
\chi_{Q[x]}=\chi_{Q[x]}\circ\chi_{S}\circ\chi_{Q[x]}=\chi_{Q[x]\Gamma S\Gamma Q[x]}(cf. \text{ Lemma } 5.12).
\end{align*}
Since $\chi_{Q[x]}(x)=1.$ Then $\chi_{Q[x]\Gamma S\Gamma Q[x]}(x)=1.$ Consequently, $x\in Q[x]\Gamma S\Gamma Q[x].$ Then there exist some $y\in S,\alpha,\beta\in\Gamma$ such that $x=x\alpha y\beta x.$ Hence $S$ is regular.
\end{proof}

\begin{theorem}
A $\Gamma$-semigroup $S$ is regular and intra-regular if and only if every quasi ideal of $S$ is idempotent.
\end{theorem}

\begin{proof}
Let $S$ be a $\Gamma$-semigroup which is regular and intra regular. Let $Q$ be a quasi ideal of $S.$ Then $Q\Gamma S\cap S\Gamma Q \subseteq  Q$. Therefore $Q\Gamma Q \subseteq Q\Gamma S\cap S\Gamma Q \subseteq Q$.
Let $a\in Q$. Then there exists $x,u,v\in S$ and $\alpha, \beta,\gamma,\delta,\sigma \in \Gamma$ such that $a=a\alpha x\beta a$ and $a=u\gamma a\delta a\sigma v$. Therefore
\begin{align*}
a&=a\alpha x\beta a=a\alpha x\beta a\alpha x\beta a=a\alpha x\beta (u\gamma a\delta a\sigma v)\alpha x\beta a\\
&=(a\alpha x\beta u\gamma a)\delta (a\sigma v\alpha x\beta a) \in Q\Gamma Q.
\end{align*}
Since $a\alpha x\beta u\gamma a\in Q\Gamma S\cap S\Gamma Q\in Q, a\sigma v\alpha x\beta a\in Q\Gamma S\cap S\Gamma Q \subseteq Q$, we deduce that $a\in Q\Gamma Q$. Therefore $Q\subseteq Q\Gamma Q$ and so $Q=Q\Gamma Q$.

Conversely, let every quasi ideal of a $\Gamma$-semigroup $S$ be idempotent. Let $a\in S.$ Now let us consider the quasi ideal $Q[a] =\{a\}\cup (a\Gamma S\cap S\Gamma a)$ generated by $a$ of $S.$ Then $Q[a]=Q[a]\Gamma Q[a]$. So $a\in Q[a]\Gamma Q[a]=(a\Gamma a)\cup (a\Gamma(a\Gamma S\cap S\Gamma a))\cup((a\Gamma S\cap S\Gamma a)\Gamma a)\cup((a\Gamma S\cap S\Gamma a)\Gamma (a\Gamma S\cap S\Gamma a)).$ Consequently $a=a\alpha x\beta a$ and $a=y\gamma a\delta a\sigma z$ for some $x,y,z\in S$ and $\alpha,\beta,\gamma,\delta,\sigma\in\Gamma$. Therefore $S$ is regular and intra regular.
\end{proof}

\begin{theorem}
Let $S$ be a $\Gamma$-semigroup. Then following are equivalent: $(1)$ Every quasi ideal of $S$ is idempotent. $(2)$ Every fuzzy quasi ideal of $S$ is idempotent.
\end{theorem}

\begin{proof}
Let $(1)$ holds. Let $\mu$ be a fuzzy quasi ideal of $S.$ Then $\mu\circ\mu\subseteq(\mu\circ\chi)\cap(\chi\circ\mu)\subseteq\mu$. Now by Theorem $5.14,$ $S$ is a regular and intra-regular $\Gamma$-semigroup. Then for $a\in S,$ $a=a\alpha x\beta a$ and $a=y\gamma a\delta a\sigma z$ for some $x,y,z\in S$ and $\alpha,\beta,\gamma,\delta,\sigma\in\Gamma$. Therefore $a=a\alpha x\beta a=a\alpha x\beta a\alpha x\beta a=a\alpha x\beta(y\gamma a\delta a\sigma z)\alpha x\beta a=(a\alpha x\beta y\gamma a)\delta (a\sigma z\alpha x\beta a).$ Now every fuzzy quasi ideal is fuzzy bi-ideal. Therefore
\begin{align*}
(\mu\circ\mu)(a)&=(\mu\circ\mu )((a\alpha x\beta y\gamma a)\delta (a\sigma z\alpha x\beta a))\\
&\geq\min\{\mu(a\alpha x\beta y\gamma a),\mu(a\sigma z\alpha x\beta a)\}\\
&\geq\min\{\min\{\mu (a),\mu (a)\},\min\{\mu (a),\mu (a)\}\}\\
&=\mu (a)\text{ for all }a\in S.
\end{align*}
Therefore $\mu\circ\mu \supseteq\mu$. Consequently $\mu =\mu\circ\mu$. Hence $(2)$ follows. Let $(2)$ holds. Let $Q$ be a quasi ideal of a $\Gamma$-semigroup $S.$ Then by Theorem $5.6,$ $C_{Q}$ is a fuzzy quasi ideal of $S.$ So by $(2),$ $C_{Q}\circ C_{Q} =C_{Q} \Rightarrow C_{Q\Gamma Q} =C_{Q}$. Therefore $Q=Q\Gamma Q.$ Hence $(1)$ follows.
\end{proof}

\begin{theorem}
Let $S$ be a $\Gamma$-semigroup. Then followings are equivalent: $(1)$ $S$ is regular and intra regular. $(2)$ Every fuzzy quasi ideal is idempotent. $(3)$ Every fuzzy bi ideal is idempotent.
\end{theorem}

\begin{proof}
$(1)\Rightarrow (3):$ Let $S$ be regular and intra-regular. Let $\mu$ be a fuzzy bi-ideal of $S.$ Let $a\in S.$ Then $a=a\alpha x\beta a$ and $a=y\gamma a\delta a\sigma z$ for some $x,y,z\in S$ and $\alpha,\beta,\gamma,\delta,\sigma\in\Gamma$. Then
\begin{align*}
(\mu\circ\mu )(a)&=(\mu\circ\mu )(a\alpha x\beta a)=(\mu\circ\mu)(a\alpha x\beta a\alpha x\beta a)\\
&=(\mu\circ\mu )(a\alpha x\beta (y\gamma a\delta a\sigma z)\alpha x\beta a)\\
&\geq\min\{\mu (a\alpha x\beta y\gamma a),\mu (a\sigma z\alpha x\beta a)\}\\
&\geq\min\{\min\{\mu(a),\mu(a)\},\min\{\mu(a),\mu(a)\}\}\\
&\geq\mu(a)\text{ for all }a\in S.
\end{align*}
Therefore $\mu\circ\mu\supseteq\mu$. Again
\begin{align*}
(\mu\circ\mu)(a)&=\underset{a=x\gamma y}{\sup}\min\{\mu(x),\mu(y)\}\\
&\leq\underset{a=x\gamma y}{\sup}\mu(x\gamma y)=\mu(a)\text{ for all }a\in S.
\end{align*}
Therefore $\mu\circ\mu\subseteq\mu$. Consequently $\mu=\mu\circ\mu$. Hence $(3)$ follows. Since every fuzzy quasi ideal is fuzzy bi-ideal, so $(3)\Rightarrow (2).$ Also by Theorem $5.14$ and $5.15$ we have $(2)\Rightarrow (1).$
\end{proof}

\begin{theorem}
Let $S$ be a $\Gamma$-semigroup and $\mu$ be a non-empty fuzzy subset of $S.$ Then $\mu$ is a fuzzy quasi ideal of $S$ if and only if $x=b\alpha s,x=t\beta c\Rightarrow\mu(x)\geq\min\{\mu(b),\mu(c)\},$ where
$b,s,t,c\in S,\alpha,\beta\in\Gamma.$
\end{theorem}

\begin{proof}
Let $\mu$ be a fuzzy quasi ideal of $S$ and $x\in S.$ Then $(\chi\circ\mu) \cap(\mu\circ\chi)\subseteq\mu$ which implies $\mu(x)\geq((\chi\circ\mu)\cap(\mu\circ\chi))(x)=\min\{(\chi\circ\mu)(x),(\mu\circ\chi)(x)\}.$ Since
$x=b\alpha s$ where $b,s\in S,\alpha\in\Gamma,$ we have
\begin{align*}
(\mu\circ\chi)(x)&=\underset{x=u\delta v}{\sup}[\min\{\mu(u),\chi(v)\}]\geq\min\{\mu(b),\chi(s)\}\\
&=\min\{\mu(b),1\}=\mu(b).
\end{align*}
Similarly, for $x=t\beta c$ where $t,c\in S,\beta\in\Gamma,$ we have $(\chi\circ\mu)(x)\geq\mu(c).$ Hence $\mu(x)\geq$min\{$\mu(b),\mu(c)\}.$

Conversely, let us suppose that $x=b\alpha s,x=t\beta c\Rightarrow\mu(x)\geq\min\{\mu(b),\mu(c)\},$ where $b,s,t,c\in S,\alpha,\beta\in\Gamma.$ If there do not exist $y,z\in S$ and $\gamma\in\Gamma$ such that $x=y\gamma z,$ then $((\chi\circ\mu)\cap(\mu\circ\chi))(x)=0\leq\mu(x).$ Then $(\chi\circ\mu)\cap(\mu\circ\chi)\subseteq\mu$ and the proof follows. Let there exist $y,m,n,z\in S$ and $\theta,\eta\in\Gamma$ such that $x=y\theta m$ and $x=n\eta z,$ Then
\begin{align*}
(\mu\circ\chi)(x)&=\underset{x=y\theta m}{\sup}[\min\{\mu(y),\chi(m)\}]\\
&=\underset{x=y\theta m}{\sup}[\min\{\mu(y),1\}]\\
&=\underset{x=y\theta m}{\sup}\mu(y).
 \end{align*}
Similarly $(\chi\circ\mu)(x)=\underset{x=n\eta z}{\sup}\mu(z).$ Then
\begin{align*}
((\chi\circ\mu)\cap(\mu\circ\chi))(x)&=\min\{(\chi\circ\mu)(x),(\mu\circ\chi)(x)\}\\
&=\min\{\underset{x=y\theta m}{\sup}\mu(y),\underset{x=n\eta z}{\sup}\mu(z)\}\\
&=\underset{x=y\theta m}{\sup}\underset{x=n\eta z}{\sup}\min\{\mu(y),\mu(z)\}\\
&\leq\underset{x=y\theta m}{\sup}\underset{x=n\eta z}{\sup}\mu(x)=\mu(x).
\end{align*}
Consequently, $(\chi\circ\mu)\cap(\mu\circ\chi)\subseteq\mu.$ Hence $\mu$ is a fuzzy quasi ideal of $S.$
\end{proof}

\begin{theorem}
Let $S$ be a $\Gamma$-semigroup and $\mu$ be a non-empty fuzzy subset of $S.$ Then $\mu$ is a fuzzy quasi ideal of $S$ if and only if $x=b\alpha s,x=t\beta c\Rightarrow\mu(x)\geq\max\{\min\{\mu(b),\mu(c)\},\min\{\mu(t),\mu(s)\}\}$ where
$b,s,t,c\in S,\alpha,\beta\in\Gamma.$
\end{theorem}

\begin{proof}
Let $\mu$ be a fuzzy quasi ideal of $S.$ Let $b,s,t,c\in S,\alpha,\beta\in\Gamma$ be such that $x=b\alpha s,x=t\beta c.$ Then by Theorem $5.17,\mu(x)\geq\min\{\mu(b),\mu(c)\}.$ Again if $x=t\beta c,x=b\alpha s,$ then by Theorem $5.17,\mu(x)\geq\min\{\mu(t),\mu(s)\}.$ Hence
\begin{align*}
\mu(x)\geq\max\{\min\{\mu(b),\mu(c)\},\min\{\mu(t),\mu(s)\}\}.
\end{align*}
Conversely, using similar argument as in Theorem $5.17$ the proof follows.
\end{proof}

\end{document}